\documentclass[10pt]{article}
\usepackage{amsmath}
\usepackage{graphicx}
\usepackage{amsfonts}
\usepackage{amssymb}

\newcommand{\1}{{\bf 1}}

\newcommand{\cb}{{\cal B}}
\newcommand{\cc}{{\cal C}}

\newcommand{\cf}{{\cal F}}

\newcommand{\cj}{{\cal J}}

\newcommand{\cac}{\mathcal C}
\newcommand{\der}{\delta}
\newcommand{\cz}{\mathcal Z}
\newcommand{\ott}{[0,T]}
\newcommand{\ist}{\int_{s}^{t}}
\newcommand{\norm}[1]{\lVert #1\rVert}
\newcommand{\ka}{\kappa}
\newcommand{\id}{\mbox{Id}}

\newcommand{\ga}{\gamma}

\newcommand{\si}{\sigma}

\newcommand{\laa}{\Lambda}

\newcommand{\R}{{\mathbb R}}

\newcommand{\G}{{\mathbb P}}

\newcommand{\lcl}{\left\{}
\newcommand{\rcl}{\right\}}
\newcommand{\lp}{\left(}
\newcommand{\rp}{\right)}
\newcommand{\lc}{\left[}
\newcommand{\rc}{\right]}

\newtheorem{theorem}{Theorem}[section]

\newtheorem{definition}[theorem]{Definition}

\newtheorem{lemma}[theorem]{Lemma}

\newtheorem{proposition}[theorem]{Proposition}

\newtheorem{remark}[theorem]{Remark}

\begin{document}
\thispagestyle{empty}

\begin{center}
\huge {

Stochastic differential 
equations with non-negativity constraints
driven by fractional Brownian motion
%with Hurst parameter H $>$ 1/2
}

\vspace{.5cm}

\normalsize {\bf Marco Ferrante$^{1}$} and {\bf Carles
Rovira$^{2,*}$}

{\footnotesize \it $^1$ Dipartimento di Matematica,
Universit\`a degli Studi di Padova,
Via Trieste 63, 35121-Padova, Italy.

$^2$ Facultat de Matem\`atiques, Universitat de Barcelona, Gran
Via 585, 08007-Barcelona.

 {\it E-mail addresses}: ferrante@math.unipd.it,
Carles.Rovira@ub.edu}

{$^{*}$corresponding author}

\end{center}

\begin{abstract}
In this paper we consider stochastic differential equations
with non-negativity constraints, driven by a fractional Brownian motion
with Hurst parameter $H>\1/2$.
We first study an ordinary integral
equation, where the integral is defined in the Young sense, and 
we prove an existence result and the boundedness of the solutions.
Then we apply this result pathwise to solve the stochastic problem.
\end{abstract}

%\vspace{0.1cm}


{\bf AMS 2000 MSC:} 60H05,
% 60H07
60H20

{\bf Running head:} stochastic differential equations with constraints

\renewcommand{\theequation}{1.\arabic{equation}}
\setcounter{equation}{0}
\section{Introduction}

The study of differential equations driven by a fractional
Brownian motion has been developed in recent years.
It has been done, using either the
formalism of the rough path analysis \cite{CQ,LQ,Gu} or the fractional calculus \cite{N-R,Z}.
As is natural, afterwards has been considered some
of the possible generalizations of the diffusion processes.
For instance, in the literature we can find now
papers about PDEs \cite{QT,CQT, DGT, GT}, Volterra equations
\cite{DT1, DT2, BR2} or systems with delay \cite{F-R, FR, NNT, LT}.

Since in some applications the quantities of interest
are naturally positive, it is also natural to consider
equations with positivity constraints.
As far as the authors know, it has only been studied up to now the case
of delay equations with positivity constraints \cite{BR}.
As we shall see, the present paper follows these steps and we shall deal with
stochastic equations with positivity constraints driven by a 
fractional Brownian motion with Hurst parameter H $>$ 1/2.
More precisely, we deal with a stochastic differential equation
with normal reflection on $\R^d$ of the form:
\begin{equation}
X(t)=x(0)+\int_0^t b(s,X(s))ds+\int_0^t \si(s,X(s))dW_s^H+Y(t),
\quad t\in(0,T],
\label{equation}
\end{equation}
where $x^i(0)>0$, for $i=1,\ldots, d$,
% where $x^i(0)\in\R^d_+=\{(x^1,\ldots,x^d)\in \R^d: x^i>0, i=1,\ldots, d\}$,
$W^H=\left\{W^{H,j},\,j=1,\ldots,m\right\}$ are independent fractional Brownian motions with Hurst parameter $H>\frac{1}{2}$ defined in a complete probability space $(\Omega,\cf,\G)$, while $Y$, the so-called regulator term,
is a vector-valued non-decreasing process
which ensures that the non-negativity constraints on X are enforced.
This can be obtained in the following customary way:
\vskip 5pt
\noindent
Set
\begin{equation}
Z(t)= x(0)+\int_0^t b(s,X(s))ds+\int_0^t \si(s,X(s))dW_s^H,\qquad t\in[0,T]. \label{num1}
\end{equation}
It is known (see e.g. \cite{DI, K-W}) that we have an explicit formula for the regulator term $Y$ in terms of $Z$, the so-called reflector term: for each $i=1,\ldots,d$
\begin{equation*}
Y^i(t)=\max_{s\in[0,t]} \left(Z^i(s)\right)^-,\qquad t\in[0,T].
\end{equation*}
Then  the solution of (\ref{equation}) satisfies
$$X(t)=
Z(t)+Y(t)\qquad t\in[0,T].$$

We call (\ref{equation}) a stochastic differential equation with reflection
driven by a fractional Brownian motion and
to the best of our knowledge this problem has not been
considered before in the wide
literature on stochastic differential equations.
\vskip 4pt

In order to deal with non-negative constraints we use the
Skorohod's mapping.
Set  $$\cc_+(\R_+,\R^d):=\left\{x\in\cc(\R_+,\R^d): x(0)\in\R^d_+\right\}.$$
\vskip 5pt
\noindent
Let us recall now the Skorokhod problem.

\begin{definition}
Given a path $z\in\cc_+(\R_+,\R^d)$, we say that a pair $(x,y)$ of functions in $\cc_+(\R_+,\R^d)$ solves the Skorokhod problem for $z$ with reflection if
\begin{enumerate}
\item $x(t)=z(t)+y(t)$ for all $t\geq 0$ and $x(t)\in\R_+^d$ for each $t\geq 0$,
\item for each $i=1,\ldots, d$, $y^i(0)=0$ and $y^i$ is nondecreasing,
\item for each $i=1,\ldots, d$, ${\displaystyle \int_0^t x^i(s)dy^i(s)=0}$ for all $t\geq 0$, so $y^i$ can increase only when $x^i$ is at zero.
\end{enumerate}
\end{definition}
\noindent
Then we have an explicit formula for $y$ in terms of $z$: for each $i=1,\ldots,d$
\begin{equation*}
y^i(t)=\max_{s\in[0,t]} \left(z^i(s)\right)^-.
\end{equation*}

The path $x$ is called the reflector of $z$ and the path $y$ is called the regulator of $z$. 
We will use the Skorokhod mapping
to force a continuous real-valued function to be 
non-negative by means of reflection at the origin. 
We will apply it to each path of $Z$ defined by (\ref{num1}).
Note that, since we are dealing with a multidimensional case, 
the mapping will be applied to each component.

\bigskip

% We have also to explain how we do understand the stochastic integral appearing
% in our stochastic equation.
At this point we have also to explain how the stochastic integral appearing
in our equation has to be understood.
Since the Hurst parameter $H > 1/2$, the stochastic integral is defined using a
pathwise approach.
We shall first consider a variation of the
Young integration theory \cite{Y} (called algebraic
integration, introduced in \cite{Gu}),
in order to define a deterministic integral
with respect to H\"older continuous function.
Then we will prove our results for
deterministic equations and at the end we will
easily apply them pathwise to
the fractional Brownian motion.

\bigskip

Let us say a few words about the strategy we have followed in order to prove our results. 
Existence and uniqueness results are usually proved together using a fixed point argument. 
In order to apply this type of argument, we have to be able to control the difference between two solutions
of our system, $\Vert x^1 - x^2 \Vert$, where $\Vert \cdot \Vert$ denotes a generic norm. 
Dealing with stochastic integrals with respect to fractional Brownian 
motion a well-posed norm to work with is the $\lambda$-H\"older one.
However, as can we seen in Remark \ref{rem13}, it is not possible to control
the difference between two regulator terms $y^1$ and $y^2$ using a  
$\lambda$-H\"older norm. 
So, we are not able to use a fixed point argument.
Actually, the existence result will be proved using an equicontinuous argument,
while the uniqueness is still an open problem. We are only able to prove the uniqueness result just up to the first
time the (up to then) unique solution has the first component being zero.

\bigskip

Here is how our paper is structured: in Section 2, we will
state our main results.
Then in Section 3 we shall recall the basic notions of
the algebraic integration theory, the Young integration
and the Skorohod mapping.
Section 4 will contain the study of the deterministic integral equations: 
the existence
and boundedness of the solutions.
Finally, Section 5 will be devoted to recall how
to apply the deterministic results to the stochastic case.

\renewcommand{\theequation}{2.\arabic{equation}}
\setcounter{equation}{0}
\section{Main results}

For any $0<\lambda\leq 1$, denote by $C^{\lambda}(s,t;\R^d)$ the space of $\lambda-$H\"older continuous functions,
namely the functions $f:[s,t]\rightarrow\R^d$ such that
\[\left\|f\right\|_{\lambda,[s,t]}:=\sup_{s\leq u<v\leq t}
\frac{\left|f(v)-f(u)\right|}{(v-u)^{\lambda}}<\infty \ ,\]
and
\[\left\|f\right\|_{\infty,[s,t]}:=\sup_{u\in[s,t]}\left|f(u)\right|.\]
\\

Let us consider the following assumptions on the coefficients.

\begin{enumerate}
\item[{\bf (H1)}]  $\sigma:[0,T]\times\R^d\rightarrow \R^{d\times m}$ is bounded and \ there exists a
constant  $K_0$ such that the
following properties hold:%
\[
\left\{
\begin{array}
[c]{l}%
i)\quad\text{Lipschitz continuity}\\
|\sigma(t,x)-\sigma(t,y)|\leq K_0|x-y|,\quad\forall\,x,y \in{}\mathbb{R}%
^{d},\,\forall\,t\in\left[  0,T\right]  \smallskip\\
% ii)\quad\text{$\nu$ H\"older continuous in time}\\
ii)\quad\text{$\nu$-H\"older continuity in time}\\
\left|  \sigma(t,x)-\sigma(s,x)\right|  \leq K_0 |t-s|^\nu, \quad\forall\,x\in
{}\mathbb{R}^{d},\,\forall\,t,s\in\left[  0,T\right].
\end{array}
\right.
\]

\item[{\bf (H2)}] $b:[0,T]\times\R^d\rightarrow \R^d$ is bounded and Lipschitz continuous in $x$, that is,  there exists a
constant  $K_0$ such that
$$
|b(t,x)-b(t,y)|\leq K_0|x-y|,\quad\forall\,x, y\in{}\mathbb{R}%
^{d},\,\forall\,t\in\left[  0,T\right] .$$

\end{enumerate}

Under these assumptions we are able to prove
that our problem admits a solution.
Our main result reads as follows:
\vskip 6pt
\begin{theorem}
Assume that
% $b$ and $\sigma$
$\sigma$ and $b$
satisfy hypothesis
{\bfseries\upshape (H1)} and {\bfseries\upshape (H2)}, respectively,
with $\nu \ge H$.
Set $\lambda_0  \in (\frac12, H)$.
Then  equation (\ref{equation}) admits a solution
\[X\in L^0(\Omega,\cf,\G;C^{\lambda_0} (0,T;\R^d)).\]
\end{theorem}

\vskip 5pt

It can also be seen that the solution has moments of any order.

\begin{theorem}\label{t1}
Assume that
% $b$ and $\sigma$
$\sigma$ and $b$
satisfy hypothesis {\bfseries\upshape (H1)} and {\bfseries\upshape (H2)}, respectively,
with $\nu \ge H$. Set $\lambda_0  \in (\frac12, H)$.
If
$X$ is a solution of (\ref{equation}),
then
\[
E ( \Vert X \Vert^p_{\lambda_0,[0,T]} ) < \infty, \quad
\forall p \ge 1.
\]
\end{theorem}

\renewcommand{\theequation}{3.\arabic{equation}}
\setcounter{equation}{0}
\section{Preliminaries}

As mentioned in the introduction, we are concerned with stochastic integral with respect to a fractional Brownian motion with Hurst parameter $H>1/2$. In order to define the stochastic integral we will use the Young integration.
 We will follow the algebraic
approach introduced in \cite{Gu} (see also \cite{GT, DT1, DT2}).
For the sake of completeness, we will recall some basic facts
and notations from
%these
those papers.
We refer the reader to the same references for
a detailed presentation.

In addition, we will recall some known results on the
Skorohod mapping and prove an inequality that we will need
 throughout the paper.

\subsection{Increments}\label{incr}

Let us begin with  the basic algebraic structures which
will allow us to define a pathwise integral with respect to
irregular functions:
first of all,  for a real number
$T>0$, a vector space $V$ and an integer $k\ge 1$ we denote by
$\cac_k(V)$ (or by $\cac_k([0,T];V)$)
 the set of continuous functions $g : [0,T]^{k} \to V$ such
that $g_{t_1 \cdots t_{k}} = 0$
whenever $t_i = t_{i+1}$ for some $i\le k-1$.
Such a function will be called a
\emph{$(k-1)$-increment}, and we will
set $\cac_*(V)=\cup_{k\ge 1}\cac_k(V)$. An elementary operator on $\cac_k(V)$
is $\der$, defined as follows:
\begin{equation}
  \label{eq:coboundary}
\delta : \cac_k(V) \to \cac_{k+1}(V), \qquad
(\delta g)_{t_1 \cdots t_{k+1}} = \sum_{i=1}^{k+1} (-1)^{k-i}
g_{t_1  \cdots \hat t_i \cdots t_{k+1}} ,
\end{equation}
where $\hat t_i$ means that this argument is omitted.
A fundamental property of $\der$
is that
$\delta \delta = 0$.
 Set $\cz\cac_k(V) = \cac_k(V) \cap \text{Ker}\delta$
and $\cb \cac_k(V) =
\cac_k(V) \cap \text{Im}\delta$.

\vspace{0.2cm}

Note that given $g\in\cac_1(V)$ and $h\in\cac_2(V)$, for any $s,u,t\in\ott$, we have
\begin{equation}
\label{eq:simple_application}
  (\der g)_{st} = g_t - g_s,
\quad\mbox{ and }\quad
(\der h)_{sut} = h_{st}-h_{su}-h_{ut}.
\end{equation}
Furthermore, it can be checked that
$\cz \cac_{k}(V) = \cb \cac_{k}(V)$ for any $k\ge 1$.
Moreover, the following property holds:
\begin{lemma}\label{exd}
Let $k\ge 1$ and $h\in \cz\cac_{k+1}(V)$. Then there exists a (non unique)
$f\in\cac_{k}(V)$ such that $h=\der f$.
\end{lemma}

Observe that Lemma \ref{exd} yields that all the elements
$h \in\cac_2(V)$ such that $\der h= 0$ can be written as $h = \der f$
for some (non unique) $f \in \cac_1(V)$.

\vspace{0.2cm}
Basically, we will use $k$-increments with $k \le 2$.
We measure the size of these increments by H\"older norms
defined in the following way: for $0\le a_1<a_2\le T$ and
$f \in \cac_2([a_1,a_2];V)$, let
$$
\norm{f}_{\mu ,[a_1,a_2]} =
\sup_{r,t\in[a_1,a_2]}\frac{|f_{rt}|}{|t-r|^\mu},
$$
and
$$
% \quad\mbox{and}\quad
\cac_2^\mu([a_1,a_2];V)=\lcl f \in \cac_2(V);\,
\norm{f}_{\mu ,[a_1,a_2]}<\infty  \rcl.
$$
Notice that the usual H\"older spaces
$\cac_1^\mu([a_1,a_2]; V)$ will be determined
        in the following way: for a continuous function
$g\in\cac_1([a_1,a_2];V)$,
 we set
\begin{equation}\label{def:hnorm-c1}
\|g\|_{\mu,[a_1,a_2]}=\|\der g\|_{\mu ,[a_1,a_2]}.
\end{equation}
We will say that $g\in\cac_1^\mu([a_1,a_2];V)$
 if $\|g\|_{\mu ,[a_1,a_2]}$ is finite.

 For $h \in \cac_3([a_1,a_2];V)$ set now
\begin{eqnarray}
  \label{eq:normOCC2}
  \norm{h}_{\gamma,\rho,[a_1,a_2]} &=& \sup_{s,u,t\in[a_1,a_2]}
\frac{|h_{sut}|}{|u-s|^\gamma |t-u|^\rho},\\
\|h\|_{\mu,[a_1,a_2]} &= &
\inf\left \{\sum_i \|h_i\|_{\rho_i,\mu-\rho_i,[a_1,a_2]} ;\, h =
 \sum_i h_i,\, 0 < \rho_i < \mu \right\} ,\nonumber
\end{eqnarray}
where the last infimum is taken over all sequences $\{h_i \in \cac_3(V) \}$
such that $h
= \sum_i h_i$ and for all choices of the numbers $\rho_i \in (0,\mu)$.
Then  $\|\cdot\|_\mu$ is  a norm on $\cac_3([a_1,a_2];V)$,
 and we set
$$
\cac_3^\mu([a_1,a_2];V):=\lcl h\in\cac_3([a_1,a_2];V);\, \|h\|_\mu<\infty \rcl.
$$
Consider $\cac_3^{1+}([a_1,a_2];V) = \cup_{\mu > 1} \cac_3^\mu([a_1,a_2];V)$.
Notice that the same kind of norms can be considered on the
spaces $\cz \cac_3([a_1,a_2];V)$, leading to the definition of the spaces
$\cz \cac_3^\mu([a_1,a_2];V)$ and $\cz \cac_3^{1+}([a_1,a_2];V)$.

\vspace{0.2cm}

The basic point in this approach to pathwise integration of irregular
processes is that, under smoothness conditions, the operator
$\delta$ can be inverted. This inverse, called $\laa$, is
defined in the following proposition, taken from \cite{LT}
and whose proof can be found in
\cite{Gu}.
\begin{proposition}
\label{prop:Lambda}
Let $0\le a_1<a_2\le T$. Then
there exists a unique linear map $\Lambda: \cz \cac^{1+}_3([a_1,a_2];V)
\to \cac_2^{1+}([a_1,a_2];V)$ such that
$$
\delta \Lambda  = \id_{\cz \cac_3^{1+}([a_1,a_2];V)}.
%\quad \mbox{ and } \quad \quad
%\Lambda  \delta= \id_{\cac_2^{1+}(V)}.
$$
In other words, for any $h\in\cac^{1+}_3([a_1,a_2];V)$ such that $\der h=0$
there exists a unique $g=\laa(h)\in\cac_2^{1+}([a_1,a_2];V)$ such that
$\der g=h$.
Furthermore, for any $\mu > 1$,
the map $\laa$ is continuous from $\cz \cac^{\mu}_3([a_1,a_2];V)$
to $\cac_2^{\mu}([a_1,a_2];V)$ and we have
\begin{equation}\label{ineqla}
\|\Lambda h\|_{\mu,[a_1,a_2]} \le \frac{1}{2^\mu-2} \|h\|_{\mu,[a_1,a_2]} ,
\qquad h \in
\cz \cac^{\mu}_3([a_1,a_2];V).
\end{equation}
\end{proposition}

\vspace{0.2cm}

\subsection{Young integration}

We will consider now the particular case where $V=\R^n$, for
an arbitrary $n\ge 1$.

Using the tools introduced in the previous subsection,
% Section \ref{incr},
here we will present a generalized integral $\ist f_u dg_u$ for
 $f \in C^{\ka}_1([0,T];\R^{n\times d})$ and
$g \in \cac_1^{\ga}([0,T];\R^d)$. Following the notations introduced
in \cite{LT, DT2},
we will sometimes write $\cj_{st}(f\, dg)$ instead of $\ist f_u dg_u$.

\vspace{0.2cm}

Let us consider first two smooth functions $f$ and $g$
defined on $\ott$. One can
write,
\begin{equation}\label{eq:def-intg-smooth-case}
\cj_{st}(f\, dg)\equiv \ist f_u \, dg_u
= f_s (\der g)_{st} + \ist (\der f)_{su} \, dg_u
=f_s (\der g)_{st} + \cj_{st}(\der f\, dg).
\end{equation}
Let us study the term $\cj(\der f\, dg)$. It is easily seen that, for
$s,u,t\in\ott$,
$$
h_{sut}\equiv \lc \der\lp \cj(\der f\, dg) \rp \rc_{sut} = (\der f)_{su}
(\der g)_{ut}.
$$
The increment $h$ is an element of $\cac_3(\R^n)$ satisfying $\der h=0$.
 Let us estimate
now the regularity of $h$: if $f\in C^{\ka}_1([0,T];\R^{n\times d})$ and $g\in
 \cac_1^{\ga}([0,T];\R^d)$,
 from
(\ref{eq:normOCC2}), it is easily checked that $h\in\cac_3^{\ga+\ka}(\R^n)$.
Hence
$h\in\cz\cac_3^{\ga+\ka}(\R^n)$, and if $\ka+\ga>1$ (which is the case if $f$
and $g$ are regular),
Proposition \ref{prop:Lambda} implies that $\cj(\der f\, dg)$
can  be written as
$$
\cj(\der f\, dg) = \laa(h) = \laa\lp  \der f \, \der g\rp,
$$
and thus, plugging this identity into (\ref{eq:def-intg-smooth-case}), we get:
\begin{equation}\label{eq:def2-intg-smooth-case}
\cj_{st}(f\, dg)
=f_s (\der g)_{st} + \laa_{st}\lp  \der f \, \der g\rp.
\end{equation}

Let us state an extension of Theorem 2.5 of \cite{LT} where it is
extended the notion of integral
whenever $f\in C^{\ka}_1([0,T];\R^{n\times d})$ and $g\in \cac_1^{\ga}([0,T];\R^d)$.

\begin{theorem}\label{t31}
Let $f \in C_1^k([0,T]; \R^{n \times d} )$ and $g \in C_1^\gamma([0,T]; \R^{d} )$ with $k + \gamma > 1$. Set
\begin{equation}
\int_s^t f dg = f_s (\delta g)_{st} + \Lambda_{st} (\delta f \delta g ). \nonumber
\end{equation}
Then:
\begin{enumerate}
\item Whenever $f$ and $g$ are smooth functions, $\int_s^t f dg $ coincides with the usual Riemann integral.

\item $\int_s^t f dg $ coincides with the Young integral as defined in \cite{Y}.

\item For any $\beta \in [0,1)$ such that $1 < \gamma + k(1-\beta)=\mu_\beta$, the generalized integral satisfies
\begin{equation}
\vert \int_s^t f dg \vert \le \Vert f \Vert_{\infty,[s,t]} \Vert g \Vert_\gamma \vert t-s\vert^\gamma + c_{\gamma,k,\beta} \Vert f \Vert_{\infty,[s,t]}^\beta
\Vert f \Vert_{k,[s,t]}^{1-\beta} \Vert g \Vert_\gamma \vert t-s \vert^{\mu_\beta}, \label{eqfon}
\end{equation}
where  $c_{\gamma,k,\beta}=2^\beta (2^{\mu_\beta}-1)^{-1}$.

\end{enumerate}

\end{theorem}

{\bf Proof:}

The proof of the original Theorem has been presented in \cite{Gu} (see also \cite{GT}, \cite{LT}).
The first two statements of our Theorem are exactly the same that those in Theorem 2.5 in \cite{LT}, so we refer the reader to this reference for their proof.

\smallskip
The last statement is a generalization of the one presented in Theorem 2.5 of \cite{LT}, where it is only considered the case $\beta=0$.  The proof for $\beta >0$ can be obtained
easily putting together the following inequality
$$
 \Vert f \Vert_{k(1-\beta),[s,t]} \leq 2^\beta  \Vert f \Vert_{\infty,[s,t]}^\beta
\Vert f \Vert_{k,[s,t]}^{1-\beta}
$$
and the inequality given in Theorem 2.5 of \cite{LT}
$$
\vert \int_s^t f dg \vert \le \Vert f \Vert_{\infty,[s,t]} \Vert g \Vert_\gamma \vert t-s\vert^\gamma + c_{\gamma,k(1-\beta)}
\Vert f \Vert_{k(1-\beta),[s,t]} \Vert g \Vert_\gamma \vert t-s \vert^{\gamma+k(1-\beta)},
$$
where $ c_{\gamma,k(1-\beta)}=(2^{\gamma+k(1-\beta)}-1)^{-1}.$

\hfill $\Box$

\subsection{Skorohod mapping}\label{sko}

We recall here from \cite{DI} a well-known result for the Skorohod mapping.

\begin{lemma} \label{le2}
For each path $z \in \cc(\R_+,\R^d)$, there exists a unique solution $(x,y)$ to the Skorokhod problem for $z$. Thus there exists a pair of functions\\ ${\displaystyle (\phi,\varphi): \cc_+(\R_+,\R^d)\rightarrow\cc_+(\R_+,\R^{2d})}$ defined by $\left(\phi(z),\varphi(z)\right)=(x,y)$. The pair $\left(\phi,\varphi\right)$ satisfies the following:
\vskip 5pt
There exists a constant $K_l>0$ such that for any $z_1,z_2\in\cc_+(\R_+,\R^d)$ we have for each $t\geq 0$,
\begin{eqnarray*}
\left\|\phi(z_1)-\phi(z_2)\right\|_{\infty,[0,t]}&\leq& K_l\left\|z_1-z_2\right\|_{\infty,[0,t]},\\
\left\|\varphi(z_1)-\varphi(z_2)\right\|_{\infty,[0,t]}&\leq& K_l\left| z_1-z_2\right\|_{\infty,[0,t]}.
\end{eqnarray*}
\end{lemma}

In our paper we will use that the
$\lambda-$H\"older norm of the regulator term $y$ is bounded
by that of $z$, as proven in the following
easy lemma.

\begin{lemma}\label{lem12}
Consider  $z \in \cc(\R_+,\R^d)$, such that $\Vert z \Vert_{\lambda,[0,T]}<\infty$. Then
for any $0\le s \le t \le T$
$$
\Vert y \Vert_{\lambda,[s,t]} \le C_d \Vert z \Vert_{\lambda,[s,t]}.$$
\end{lemma}
{\bf Proof:}

\noindent
Take $u,v$ such that $ s \le u <v \le t$. Fixed a component $i$, we wish to study
$$ \frac{\vert y_v^i - y_u^i \vert}{(v-u)^\lambda}.$$
When $y_u^i=y_v^i$, this is clearly zero.
On the other hand, when $y_v^i > y_u^i$, let us define
\begin{eqnarray*}
u^*&:=&\sup\{ u'\ge u; y^i_{u'}=y^i_{u} \},\\
v^*&:=&\inf\{ v' \le v; y^i_{v'}=y^i_{v} \}.
\end{eqnarray*}
Then, $u \le u^* < v^* \le v$ and $y_u^i=y_{u^*}^i, y_v^i=y_{v^*}^i$. So
$$ \frac{\vert y_v^i - y_u^i \vert}{(v-u)^\lambda} \le \frac{\vert y_{v^*}^i - y_{u^*}^i \vert}{(v^*-u^*)^\lambda} =
\frac{\vert z_{v^*}^i -z_{u^*}^i \vert}{(v^*-u^*)^\lambda} $$
where the last equality follows from the fact that
$y^i$ and $z^i$ coincides whenever $y^i$ is not constant.

Then, note that
$$
\sup_{s \le u <v \le t} \frac{\vert y_v^i - y_u^i \vert}{(v-u)^\lambda}   \le
\sup_{s \le u^* <v^* \le t} \frac{\vert z_{v^*}^i -z_{u^*}^i \vert}{(v^*-u^*)^\lambda} \leq
 \Vert z \Vert_{\lambda,[s,t]}.$$
Finally, we get that
$$
\Vert y \Vert_{\lambda,[s,t]} \leq \Big( \sum_{i=1}^d \big( \sup_{s \le u <v \le t} \frac{\vert y_v^i - y_u^i \vert}{(v-u)^\lambda} \big)^2 \Big)^\frac12   \le d^\frac12
 \Vert z \Vert_{\lambda,[s,t]}.$$

\hfill $\Box$

\begin{remark}\label{rem13}
It is possible to prove that a similar estimate
does not hold for the difference of two regulator terms, 
result that we would need in order to prove a uniqueness theorem
in the H\"older norm framework.
Indeed,
let $0<t_1<t_2<t$, $\lambda \in(0,1)$, and
take $z^1, z^2\in C^\lambda([0,t])$
defined as
\[
z^1(s) = [(t_2-s)/(t_2-t_1)-1] \1_{(t_1,t_2]}(s)-\1_{(t_2,t]}(s)
\]
\[
z^2(s) = s/t_1 \1_{[0,t_1]}(s)+(t_2-s)/(t_2-t_1) \1_{(t_1,t_2]}(s)
\]
(note that $z^1(0)=z^2(0)$).
It is easy to see that
$y^1(s) = [1-(t_2-s)/(t_2-t_1)] \1_{(t_1,t_2]}(s)+\1_{(t_2,t]}(s)$,
while $y^2(s) \equiv 0 $. We get then
\[
\|y^2-y^1\|_{\lambda,[0,t]}=\|y^1\|_{\lambda,[0,t]}
= \frac{1}{(t_2-t_1)^\lambda}
\]
while
\[
\|z^2-z^1\|_{\lambda,[0,t]}=
\frac{1}{(t_1)^\lambda}
\]
Taking $t_1$ fixed and $t_2-t_1$ small, we prove that
in general
the $\lambda-$H\"older norm of the difference of two
regulator terms cannot be bounded by the $\lambda-$H\"older norm
of the difference of $z^1$ and $z^2$.
\end{remark}

\renewcommand{\theequation}{4.\arabic{equation}}
\setcounter{equation}{0}
\section{Deterministic integral equations}

In this section we will prove all the deterministic results.
\vskip 5pt
\noindent
Consider the deterministic
% stochastic
differential equation on $\R^d$
\begin{equation}
x(t)=x(0)+\int_0^t b(s,x(s))ds+\int_0^t\sigma(s,x(s))dg_s+ y(t),\quad t\in (0,T],\label{determinista}
\end{equation}
{where for each} $i=1,\ldots,d$
\begin{equation*}
{\mathit y^i(t)=\max_{s\in[0,t]} \left(z^i(s)\right)^-,\qquad t\in[0,T],}
\end{equation*}
{and}
\begin{equation*}
{\mathit z(t)= x(0)+\int_0^t b(s,x(s))ds+\int_0^t \si(s,x(s))dg_s,\qquad t\in[0,T].}
\end{equation*}

We will assume that the driving noise $g$ belongs to
$C^\gamma([0,T]; \R^m)$ with $\gamma > \frac12$.
Then, the integral with respect to $g$ has
to be interpreted in the Young
sense and we will find a solution
$x$ in the space $C^\lambda([0,T]; \R^d)$
with $\lambda \in (\frac12, \gamma)$.

The result of existence reads as follows.

\begin{theorem}\label{teoexi}
Assume that
% $b$ and $\sigma$
$\sigma$ and $b$
satisfy hypothesis {\bfseries\upshape (H1)}
and {\bfseries\upshape (H2)}, respectively,
with $\nu \ge \gamma$ .
Set $\lambda \in (\frac12,\gamma)$.
Then  equation (\ref{determinista}) has
% an unique
a solution $x\in C^{\lambda}([0,T];\R^d_+)$.
\end{theorem}
{\bf Proof:}

To prove that equation (\ref{determinista}) admits a solution on $[0,T]$, we shall prove first that it has a solution on $[0,T_1]$ for $T_1$ small enough ($T_1$ will be defined later). Then we will extend the solution to $[0,T]$ using an induction argument to extend the result from $[0,nT_1]$ to  $[0,(n+1)T_1]$.

\vskip 7pt
\noindent
{\it STEP 1}: Study on $[0,T_1]$.

Let us consider
\begin{equation}
 x^{(1)}(t)=z^{(1)}(t)=x(0); y^{(1)}(t)=0 \qquad t\in [0,T],
\end{equation}
and for all $n >1$
\begin{equation}\label{ambn}
 x^{(n+1)} (t)=x(0)+\int_0^t b(s,x^{(n)}(s))ds+\int_0^t\sigma(s,x^{(n)}(s))dg_s+y^{(n)}(t),
\qquad  t\in[0,T],
\end{equation}
 where for each $i=1,\ldots,d$
\begin{equation*}
 y^{(n),i} (t)=\max_{s\in[0,t]} \left(z^{(n),i}(s)\right)^-,\qquad t\in[0,T],
\end{equation*}
with
\begin{equation*}
 z^{(n)}(t)= x(0)+\int_0^t b(s,x^{(n)}(s))ds+\int_0^t \si(s,x^{(n)}(s))dg_s,\qquad t\in[0,T].
\end{equation*}
\vskip 2pt

\noindent {\it Step 1.1:  Properties of the functions $x^{(n)}$}

It follows that $x^{(n)} \in C^\lambda([0,T_1],\R_+^d)$ for all $n \ge 1$.
Indeed, from Lemma $\ref{lem12}$,
we have that
\begin{equation*}
\Vert x^{(n+1)} \Vert_{\lambda,[0,T_1]} \le \Vert z^{(n)} \Vert_{\lambda,[0,T_1]} + \Vert y^{(n)} \Vert_{\lambda,[0,T_1]}
\le C_d \Vert z^{(n)} \Vert_{\lambda,[0,T_1]}.
\end{equation*}
Using Theorem \ref{t31} and the hypothesis on the coefficients
\begin{eqnarray*}
 & & \vert \int_s^t b(u,x^{(n)}(u))du+\int_s^t \si(u,x^{(n)}(u))dg_u \vert \\
 & & \qquad \le \Vert b \Vert_\infty \vert t-s \vert + \Vert \sigma \Vert_\infty \Vert g \Vert_\gamma \vert t-s \vert^\gamma  \\
 & & \qquad \qquad \qquad +
 c_{\gamma,\lambda}  \Vert g \Vert_\gamma \Vert \sigma(.,x^{(n)}(.)) \Vert_{\lambda,[s,t]} \vert t-s \vert^{\gamma+\lambda}.
 \end{eqnarray*}
 Using that
 $$
 \Vert \sigma(.,x^{(n)}(.)) \Vert_{\lambda,[s,t]} \le K_0 \left( \Vert x^{(n)} \Vert_{\lambda,[s,t]} + \vert t-s \vert^{\nu-\lambda} \right),$$
 we can write that
 \begin{equation}\label{desx}
 \Vert x^{(n+1)} \Vert_{\lambda,[s,t]} \le  h(t-s) + M_1 \Vert x^{(n)} \Vert_{\lambda,[s,t]} \vert t-s \vert^\gamma,
 \end{equation}
 where $h(t)=C_d(\Vert b \Vert_\infty t^{1-\lambda}  + \Vert \sigma \Vert_\infty \Vert g \Vert_\gamma t^{\gamma-\lambda} +  c_{\gamma,\lambda} K_0 \Vert g \Vert_\gamma t^{\gamma-\lambda+\nu}), \break
 M_1= C_d c_{\gamma,\lambda} K_0 \Vert g \Vert_\gamma.$
 Repeating iterativily inequality (\ref{desx}) with $s=0$ we get that
 \begin{equation}\label{desx2}
 \Vert x^{(n+1)} \Vert_{\lambda,[0,t]} \le  h(t) \left( \sum_{i=0}^{n-1} ( M_1 t^\gamma )^{i} +  ( M_1 t^\gamma )^{n}
 \Vert x^{(1)} \Vert_{\lambda,[0,t]} \right)=   h(t)  \sum_{i=0}^{n-1} ( M_1 t^\gamma )^{i} .
 \end{equation}
So, choosing $T_1$ such that
$T_1 < \left(1/M_1\right)^{1/\gamma}$,
\begin{equation}\label{fitaholder}
\sup_{n \ge 1} \Vert x^{(n)} \Vert_{\lambda,[0,T_1]} := K_1 < \infty.
\end{equation}
Since $x^{(n)}(0)=x(0)$ for all $n$, it follows easily that
\begin{equation}\label{fitainf}
\sup_{n \ge 1} \Vert x^{(n)} \Vert_{\infty,[0,T_1]} := K_2 < \infty.
\end{equation}

\noindent  {\it Step 1.2: Definition of the solution}

The sequence of functions $x^{(n)}$ is equicontinuous and bounded in $C([0,T_1];\R^d)$. Therefore there exists a subsequence  $\{x^{(n_j)}\}_{j\ge 1}$ that converges uniformly to a
function $x \in  C([0,T_1];\R^d)$ .

Moreover $x$ belongs to $C^\lambda([0,T_1],\R_+^d)$. Indeed, fixed $\varepsilon >0$ let us choose $n_j $ such that
$\Vert x-x^{n_j} \Vert_{\infty,[0,T_1]} \le \varepsilon$. Then, for all $s,t \in [0,T_1]$
\begin{eqnarray*}
\vert x (t)- x(s) \vert & \le & \vert x(t) - x^{(n_j)}(t) \vert +  \vert x^{(n_j)}(t) - x^{(n_j)}(s) \vert +  \vert x^{(n_j)}(s) - x (s) \vert \\
&\le & 2\varepsilon + K_1 \vert t-s\vert^\lambda.
\end{eqnarray*}
Since this inequality is true for all $\varepsilon >0$,
% letting $\varepsilon$ to $0$,
we get that $x \in  C^\lambda ([0,T_1];\R_+^d)$ and, for any $ t \in [0,T_1]$,  the Young integral
$$\int_0^t \sigma(s,x(s)) dg_s $$
is well-defined.

\vskip 2pt

\noindent  {\it Step 1.3: $x$ satisfies (\ref{determinista}) on $[0,T_1]$.}

Since  $\{x^{(n_j)}\}_{j\ge 1}$  converges uniformly to $x$  and $b$ is Lipschitz in space, we get that
\begin{equation}\label{limles}
\lim_{n_j \to \infty} \Vert \int_0^. b(s,x^{(n_j)}(s)) ds - \int_0^. b(s,x(s)) ds \Vert_{\infty,[0,T_1]} = 0.
\end{equation}

On the other hand, using the hypothesis on the coefficients and Theorem \ref{t31}, for any $t \in [0,T_1]$:
\begin{eqnarray*}
& &
\left|\int_0^t \left(\sigma(s,x(s))-\sigma(s,x^{(n_j)}(s))\right)dg_s\right|
\\&\leq&
K_0 \Vert x^{(n)}-x^{(n-1)}\Vert_{\infty,[0,T_1]}  \Vert g \Vert_{\gamma} T_1^\gamma\\
& & \qquad +  K_0^\beta c_{\gamma,\lambda,\beta}  \Vert x-x^{(n_j)}\Vert_{\infty,[0,T_1]}^\beta
 \Vert \sigma(.,x(.))-\sigma(.,x^{(n_j)}(.)) \Vert_{\lambda,[0,T_1]}^{1-\beta} \Vert g \Vert_{\gamma} T_1^{\mu_\beta}.
\end{eqnarray*}
Since $ \Vert \sigma(.,x(.))-\sigma(.,x^{(n_j)}(.)) \Vert_{\lambda,[0,T_1]}<\infty,$
\begin{equation}\label{limyou}
\lim_{n_j \to \infty} \Vert \int_0^.\sigma(s,x(s)) dg_s - \int_0^. \sigma(s,x^{(n_j)}(s)) dg_s \Vert_{\infty,[0,T_1]} = 0.
\end{equation}

Finally, for each $i=1,\ldots,d$, set
\begin{equation*}
 y^{i} (t)=\max_{s\in[0,t]} \left(z^{i}(s)\right)^-,\qquad t\in[0,T_1],
\end{equation*}
where
\begin{equation*}
 z(t)= x_0+\int_0^t b(s,x(s))ds+\int_0^t \si(s,x(s))dg_s,\qquad t\in[0,T_1].
\end{equation*}
Using Lemma \ref{le2} we have
\[
\sup_{t\in[0,T_1]} \left|y(t)-y^{(n_j)}(t)\right|\leq K_l \sup_{t\in[0,T_1]} \left|z(t)-z^{(n_j)}(t)\right|.\] From (\ref{limles}) and
(\ref{limyou}) we obtain now that
\begin{equation}\label{limy}
\lim_{n_j \to \infty} \Vert y- y^{(n_j)}  \Vert_{\infty,[0,T_1]} = 0.
\end{equation}

Letting $n_j$ to infinity in (\ref{ambn}) and using  (\ref{limles}) ,(\ref{limyou}) and (\ref{limy}), we get that $x$ satisfies (\ref{determinista}).

\vskip 7pt
\noindent
{\it STEP 2.}:  We will assume that we have defined the solution on $[0,NT_1]$. We will show first the extension to $[NT_1,(N+1)T_1]$  (assuming $(N+1)T_1 \le T$).

Let $x$ be a solution defined in $[0,N T_1]$. Then, let
 us define
 \begin{eqnarray*}
 && x^{(1)} (t) := x(t) \1_{[0,NT_1]}(t) + x(T_1) \1_{(NT_1,T]}(t),\\
&& z^{(1)} (t) := z(t) \1_{[0,NT_1]}(t) + z(T_1) \1_{(NT_1,T]}(t).
\end{eqnarray*}
Moreover, for all $n >1$
\begin{eqnarray*}%\label{ambnmes}
&& x^{(n+1)} (t)= x(t), \qquad\quad t \in [0,NT_1],\\
&& x^{(n+1)} (t)=z(NT_1) +\int_{NT_1}^t b(s,x^{(n)}(s))ds+\int_{N T_1}^t\sigma(s,x^{(n)}(s))dg_s+y^{(n)}(t),\\
&& \qquad\qquad\qquad\qquad\qquad\qquad\qquad\qquad\qquad\qquad \quad t \in (NT_1,T],
\end{eqnarray*}
 where for each $i=1,\ldots,d$
\begin{equation*}
 y^{(n),i} (t)=\max_{s\in[0,t]} \left(z^{(n),i}(s)\right)^-,\qquad t\in[NT_1,T],
\end{equation*}
with
\begin{eqnarray*}
&& z^{(n)}(t)= z(t), \qquad\quad t \in [0,T_1],\\
&& z^{(n)}(t)= z(NT_1)+\int_{NT_1}^t b(s,x^{(n)}(s))ds+\int_{NT_1}^t \si(s,x^{(n)}(s))dg_s,,\\
&& \qquad\qquad\qquad\qquad\qquad\qquad\qquad\qquad\qquad\qquad \quad t \in (NT_1,T].
\end{eqnarray*}
Note that
$$ y^{(1)} (t) := y(t) \1_{[0,NT_1]}(t) + y(T_1) \1_{(NT_1,T]}(t).$$

Repeating the same computations given in {\it Step 1.1} we get that
 \begin{equation}
 \Vert x^{(n+1)} \Vert_{\lambda,[N T_1,t]} \le  h(t-N T_1)  \sum_{i=0}^{n-1} ( M_1 (t-N T_1)^\gamma )^{i} ,
 \end{equation}
 where $h(t)$ and $M_1$ are the same that appear in (\ref{desx2}). Using the same ideas employed in {\it Step 1.1}, we obtain that
\begin{equation*}
\sup_{n \ge 1} \Vert x^{(n)} \Vert_{\lambda,[NT_1,(N+1)T_1]} := K_1 < \infty,
\end{equation*}
and
\begin{equation*}
\sup_{n \ge 1} \Vert x^{(n)} \Vert_{\infty,[NT_1,(N+1)T_1]} := K_2 < \infty.
\end{equation*}

Following now the method used in {\it Step 1.2} and {\it Step 1.3},  there exists a subsequence  $\{x^{(n_j)}\}_{j\ge 1}$ that converges uniformly to a
function $x^{[2]} \in  C([NT_1,(N+1)T_1];\R^d)$ . Moreover $x^{[2]}$ belongs to $C^\lambda([NT_1,(N+1)T_1],\R_+^d)$ and $x^{[2]}$ satisfies (\ref{determinista}) in $[NT_1,(N+1)T_1]$.

Set now:
$$x(t)= x(t) \1_{[0,NT_1]} (t) +  x^{[2]}(t) \1_{(NT_1,(N+1)T_1]} (t).$$
Clearly,  $x$ belongs to $C^\lambda([0,(N+1)T_1],\R^d)$ and $x$ satisfies (\ref{determinista}) in $[0,(N+1)T_1]$.

\hfill $\Box$

\vskip 5pt

\begin{remark} 
The study of the uniqueness of the solution is an open problem, due to the
fact that we are not able to bound the H\"older norm of the difference
of regulator terms with the same norm of the difference of the reflected terms.
We can only get the uniqueness in a small time interval and 
when $\sigma$ 
does not depend on time. Assuming that
$\sigma:\R\rightarrow \R^{ m}$ is bounded and Lipschitz
continuous,  $b$
satisfy hypothesis {\bfseries\upshape (H2)},  $\nu \ge \gamma$ and
$\lambda_0  \in (\frac12, \gamma)$, 
then there exists 
$\varepsilon>0$ such that equation (\ref{determinista}) 
has a unique solution on $[0,\varepsilon]$
and
$x\in C^{\lambda_0}([0,\varepsilon];\R^d_+)$.

Indeed, fix $x$ a solution of our equation in $[0,T]$.
Since $x^i(0)>0$ for any $i$, there exists $\varepsilon>0$
such that $x^i(t)>0$ for any $i$ and $t\in[0,\varepsilon]$.
Let $x_1$ be a second solution of our equation and let 
$\delta:=\sup\{t\ge0:x_1^i(t)>0, i=1,\ldots, k\}$. 
Clearly $\delta>0$ and assume, w.l.g. that $\delta\le \varepsilon$.
On $[0,\delta)$, both solutions $x$ and $x_1$ are nonnegative and $y(s)=y_1(s)=0$
for all $s \in  [0,\delta)$. For any $t_0 \in  [0,\delta)$
\begin{eqnarray}
\nonumber
\Vert x - x _1 \Vert_{\lambda,[0,t_0]} & \leq &
\Vert \int_0^.  (b(u,x^a(u)) - b(u,x^b(u))) du  \Vert_{\lambda,[0,t_0]}\\
& & +
\Vert \int_0^.  (\sigma(x(u)) - \sigma(x_1(u))) dg_u  \Vert_{\lambda,[0,t_0]}\label{es1}.
\end{eqnarray}
On one hand
\begin{eqnarray}\nonumber
& &\Vert \int_0^.  (b(u,x(u)) - b(u,x_1(u))) du  \Vert_{\lambda,[0,t_0]} \\ \nonumber  &  & \qquad
%  \leq
=
\sup_{0 \le s<t \le t_0} \frac{\vert \int_s^t  (b(u,x(u)) - b(u,x_1(u)))
du \vert}{\vert t-s \vert^\lambda}\\
&  & \qquad \leq  K_0 \Vert x - x_1 \Vert_{\infty,[0,t_0]}  t_0^{1-\lambda}
\leq  K_0 \Vert x - x_1 \Vert_{\lambda,[0,t_0]}  t_0\label{es2}.
\end{eqnarray}
On the other hand, using Theorem \ref{t31}
\begin{eqnarray}\nonumber
& &\Vert \int_0^.  (\sigma(x(u)) - \sigma(x_1(u))) dg_u   
\Vert_{\lambda,[0,t_0]} \\ \nonumber &  & \qquad \leq
 K_0 \Vert x- x_1 \Vert_{\infty,[0,t_0]} \Vert g \Vert_\gamma
 t_0^{\gamma-\lambda} + c_{\gamma,\lambda}
\Vert x- x_1 \Vert_{\lambda,[0,t_0]} \Vert g \Vert_\gamma t_0^{\gamma}
\\  &  & \qquad \leq
 ( K_0 + c_{\gamma,\lambda})
\Vert x- x_1 \Vert_{\lambda,[0,t_0]} \Vert g \Vert_\gamma t_0^{\gamma}
\label{es3}.
\end{eqnarray}
Putting (\ref{es2}) and   (\ref{es3})  in  (\ref{es1}) we get that
$$
\Vert x^a- x^b \Vert_{\lambda,[0,t_0]} \leq
\Big( K_0 t_0 +  ( K_0 + c_{\gamma,\lambda})
\Vert g \Vert_\gamma t_0^{\gamma} \Big)
\Vert x- x_1 \Vert_{\lambda,[0,t_0]} .$$
Choosing $t_0$ small enough it follows that $\Vert x- x_1 \Vert_{\lambda,[0,t_0]}=0$. 
Since $x(0)=x_1(0)$ it follows that
$x=x_1$ on $[0,t_0]$.
Since we can repeat this arguments on $[t_0, 2t_0]$
and  further on, we get that
$x=x_1$ on $[0,\delta]$.

\ From the continuity of our solutions we have that $x(\delta) = x_1(\delta)=0$,
 $y(\delta) = y_1(\delta)=0$ and  $\delta=\varepsilon$.

\end{remark}

\bigskip

We will finish the study of the deterministic case obtaining a bound for the H\"older norm of the solutions.

\bigskip

\begin{theorem}\label{teofit}
Assume that
% $b$ and $\sigma$
$\sigma$ and $b$ satisfy hypothesis {\bfseries\upshape (H1)}
and {\bfseries\upshape (H2)}, respectively,
with $\nu \ge \gamma$ and
set $\lambda \in (\frac12,\gamma)$.
Given $x$ a solution of equation (\ref{determinista}), it holds that
$\Vert  x \Vert_{\lambda,[0,T]} \le M_2 + M_3 \Vert g \Vert_\gamma^\frac{1}{\gamma}$,
where $M_2$ and $M_3$ are positive constants not depending on $g$.
\end{theorem}

{\bf Proof:}

Using that $b$ and $\sigma$ are bounded, Lemma \ref{lem12},  Theorem \ref{t31} and that
$$
\Vert \sigma(.,x(.)) \Vert_{\lambda,[s,t]}\le K_0 \vert t-s \vert^{\nu-\lambda} + K_0 \Vert x
 \Vert_{\lambda,[s,t]},$$
we get that for any $s \le t$
\begin{eqnarray*}
\Vert x
 \Vert_{\lambda,[s,t]} & \le &  \Vert z
 \Vert_{\lambda,[s,t]} +  \Vert y
 \Vert_{\lambda,[s,t]} \le  (C_d +1) \Vert z
 \Vert_{\lambda,[s,t]}\\
 & \le &
(C_d +1) \left( \Vert b \Vert_\infty \vert t-s \vert^{1-\lambda}
+ \Vert \sigma \Vert_\infty \Vert g \Vert_\gamma \vert t-s\vert^{\gamma-\lambda} \right.
\\
& &  \qquad \left.  + c_{\gamma,\lambda}
\Vert \sigma(.,x(.)) \Vert_{\lambda,[s,t]} \Vert g \Vert_\gamma \vert t-s \vert^{\gamma}\right) \\
& \le &
(C_d +1) \left( \Vert b \Vert_\infty \vert t-s \vert^{1-\lambda}
+ (\Vert \sigma \Vert_\infty + T^\nu K_0 c_{\gamma,\lambda}) \Vert g \Vert_\gamma \vert t-s\vert^{\gamma-\lambda} \right.
\\
& &  \qquad \left.  + K_0 c_{\gamma,\lambda}
\Vert x \Vert_{\lambda,[s,t]} \Vert g \Vert_\gamma \vert t-s \vert^{\gamma}\right).
\end{eqnarray*}

Set $M_{d,\gamma,\lambda}:= (C_d+1) K_0 c_{\gamma,\lambda} $. So, for any $0 \le s < t \le T$ such that
\begin{equation} \label{dests} M_{d,\gamma,\lambda}
\Vert g \Vert_\gamma \vert t-s \vert^{\gamma} \le \frac12,
\end{equation}
we get that
\begin{equation}\label{petita}
\Vert x
 \Vert_{\lambda,[s,t]} \le 2 (C_d +1) \left( \Vert b \Vert_\infty \vert t-s \vert^{1-\lambda}
+ (\Vert \sigma \Vert_\infty + T^\nu K_0 c_{\gamma,\lambda}) \Vert g \Vert_\gamma \vert t-s\vert^{\gamma-\lambda} \right)
\end{equation}

%More precisely, we can also obtain the following bound:
%\begin{eqnarray}
%\Vert x
% \Vert_{\lambda,[s,t]} &\le& 2 (C_d +1)  \left( \Vert b \Vert_\infty \vert t-s \vert^{1-\lambda}
%+ (\Vert \sigma \Vert_\infty +  T M_1 c_{\gamma,\lambda}) \Vert g \Vert_\gamma \left( 2 M_{d,\gamma,\lambda}
%\Vert g \Vert_\gamma \right)^{-\frac{1}{\gamma}} )^{\gamma-\lambda} \right)\nonumber\\
%&\le& 2 (C_d +1)  \left( \Vert b \Vert_\infty \vert t-s \vert^{1-\lambda}
%+ (\Vert \sigma \Vert_\infty +  T M_1 c_{\gamma,\lambda})  \left( 2 M_{d,\gamma,\lambda}
%\right)^{\frac{\lambda}{\gamma}-1 } \Vert g \Vert_\gamma^{\frac{\lambda}{\gamma}}\right)\label{primera}
%\end{eqnarray}

Note that given any $0 \le s < t \le T$ that do not satisfy (\ref{dests}), we can choose $t_0=s < t_1 < \ldots  < t_n =t$
such that for all $ i \in {1,\ldots,n}$
$$ M_{d,\gamma,\lambda}
\Vert g \Vert_\gamma \vert t_{i}-t_{i-1} \vert^{\gamma} \le \frac12$$
with
$$ n \le 2 \vert t-s \vert \left( 2 M_{d,\gamma,\lambda}
\Vert g \Vert_\gamma \right)^{\frac{1}{\gamma}}.$$

Then, using (\ref{petita}), we have that
\begin{eqnarray*}
& &\frac{\vert x(t) - x(s) \vert}{\vert t-s \vert^\lambda}
\le
\sum_{i=1}^n \frac{\vert x(t_{i}) - x(t_{i-1}) \vert}{\vert t_{i}-t_{i-1} \vert^\lambda}
\frac{\vert t_{i}-t_{i-1} \vert^\lambda}{\vert t-s \vert^\lambda}\\
& & \qquad \le \sum_{i=1}^n \Vert x
 \Vert_{\lambda,[t_{i-1},t_i]} \frac{\vert t_{i}-t_{i-1} \vert^\lambda}{\vert t-s \vert^\lambda} \\
& & \qquad \le \sum_{i=1}^n \frac{2 (C_d +1) \left( \Vert b \Vert_\infty \vert t_i-t_{i-1} \vert
+ (\Vert \sigma \Vert_\infty + T^\nu K_0 c_{\gamma,\lambda}) \Vert g \Vert_\gamma \vert t_i-t_{i-1}\vert^{\gamma} \right)}{\vert t-s \vert^\lambda}\\
 & & \qquad \le 2 (C_d +1)  \Vert b \Vert_\infty  T^{1-\lambda} \\
 & & \qquad \quad +
 \sum_{i=1}^n \frac{2 (C_d +1) \left((\Vert \sigma \Vert_\infty + T^\nu K_0 c_{\gamma,\lambda}) \Vert g \Vert_\gamma \vert t_i-t_{i-1}\vert^{\gamma} \right)}{\vert t-s \vert^\lambda}\\
 & & \qquad \le 2 (C_d +1)  \Vert b \Vert_\infty  T^{1-\lambda} \\
 & & \qquad \quad +
 \left( 2 \vert t-s \vert \left( 2 M_{d,\gamma,\lambda}
\Vert g \Vert_\gamma \right)^{\frac{1}{\gamma}}\right)   \frac{2 (C_d +1) \left((\Vert \sigma \Vert_\infty + T^\nu K_0 c_{\gamma,\lambda}) \Vert g \Vert_\gamma \right)}{\vert t-s \vert^\lambda  (2 M_{d,\gamma,\lambda} \Vert g \Vert_\gamma) }\\
& & \qquad \le 2 (C_d +1)  \Vert b \Vert_\infty  T^{1-\lambda} \\
 & & \qquad \quad +
 T^{1-\lambda}  2^{1+\frac{1}{\gamma}}M_{d,\gamma,\lambda}^{\frac{1}{\gamma}-1}  (C_d +1) \left(\Vert \sigma \Vert_\infty + T^\nu K_0 c_{\gamma,\lambda}  \right) \Vert g \Vert_\gamma^\frac{1}{\gamma}.
\end{eqnarray*}

\ From this last inequality, it follows easily that
$\Vert  x \Vert_{\lambda,[0,T]} \le M_2 + M_3 \Vert g \Vert_\gamma^\frac{1}{\gamma}$.

\hfill $\Box$

\section{ Stochastic integral equations}

In this section we
 apply the deterministic results
 to prove the main theorems of this paper.

The following Lemma, taken from \cite{N-R} (see Lemma 7.4) is basic in order to extend the deterministic results to the stochastic ones.

\begin{lemma}\label{f1}
Let $\{W_t^H; t \ge 0 \}$ be a fractional Brownian
motion of Hurst parameter $H \in (0,1).$
Then for
every $0 < \varepsilon < H$ and $T>0$ there exists a positive random variable $\eta_{\varepsilon,T}$ such that
$E (\vert \eta_{\varepsilon,T} \vert^p) < \infty$ for all $p \in [1,\infty)$ and for all $s,t \in [0,T]$
$$
\vert W_t^H - W_s^H \vert \le \eta_{\varepsilon,T} \vert t-s \vert^{H-\varepsilon},\qquad {\rm a.s}.$$
\end{lemma}

The stochastic integral appearing throughout this paper $ \int_0^T
u(s) dW^H_s $ is a Young integral. This integral exists if the process $u(s)$ has
H\"older continuous trajectories of order larger than $\lambda$ such that $H + \lambda > 1$.

On the other hand, notice that from lemma \ref{f1}, for any $\gamma < H$ it holds that
$$
E( \Vert W^H \Vert_\gamma^p ) < \infty,$$
for all $p \in [1,\infty)$.

\bigskip

Choosing $\gamma$ and $\lambda$, such that $\gamma=H - \varepsilon_1$ with $\gamma > \frac12$ and $\lambda \in (\frac12, \gamma)$,
Theorems \ref{teoexi} and \ref{teofit}
follows almost surely.
Our stochastic theorems of existence and boundedness 
of the moments hold then clearly.

\section*{Acknowledgements}

We would like to thank A. Deya for his comments that have helped us to simplify the proof of Theorem \ref{t31}.

This work was partially supported by
COFIN 2008-YYYBE4-002 of MIUR (Marco Ferrante)
and
DGES Grant MTM09-07203 (Carles Rovira).

\end{document}